\documentclass[12pt]{amsart}
\usepackage{amscd}
\def\NZQ{\Bbb}               
\def\NN{{\NZQ N}}

\def\ZZ{{\NZQ Z}}

\def\frk{\frak}               

\def\pp{{\frk p}}

\def\mm{{\frk m}}

\def\Phi{{\frk n}}
\def\Phi{{\frk N}}
\def\opn#1#2{\def#1{\operatorname{#2}}} 
\opn\chara{char} \opn\length{\ell} \opn\pd{pd} \opn\rk{rk} \opn\projdim{proj\,dim} \opn\injdim{inj\,dim} \opn\rank{rank} \opn\depth{depth} \opn\grade{grade}
\opn\height{height} \opn\embdim{emb\,dim} \opn\codim{codim}

\opn\Tr{Tr} \opn\bigrank{big\,rank} \opn\superheight{superheight}\opn\lcm{lcm}
\opn\trdeg{tr\,deg}
\opn\reg{reg} \opn\lreg{lreg} \opn\ini{in} \opn\lpd{lpd} \opn\size{size}\opn{\mult}{mult}
\opn\div{div} \opn\Div{Div} \opn\cl{cl} \opn\Cl{Cl}
\opn\Spec{Spec} \opn\Supp{Supp} \opn\supp{supp} \opn\Sing{Sing} \opn\Ass{Ass} \opn\Min{Min}
\opn\Ann{Ann} \opn\Rad{Rad} \opn\Soc{Soc}
\opn\Syz{Syz} \opn\Im{Im} \opn\Ker{Ker} \opn\Coker{Coker} \opn\Am{Am} \opn\Hom{Hom} \opn\Tor{Tor} \opn\Ext{Ext} \opn\End{End} \opn\Aut{Aut} \opn\id{id}

\opn\nat{nat}
\opn\pff{pf}
\opn\Pf{Pf} \opn\GL{GL} \opn\SL{SL} \opn\mod{mod} \opn\ord{ord} \opn\Gin{Gin} \opn\Hilb{Hilb}\opn\adeg{adeg}\opn\std{std}\opn\ip{infpt}
\opn\aff{aff} \opn\con{conv} \opn\relint{relint} \opn\st{st} \opn\lk{lk} \opn\cn{cn} \opn\core{core} \opn\vol{vol} \opn\link{link} \opn\star{star}
\opn\gr{gr}

\def\pot#1#2{#1[\kern-0.28ex[#2]\kern-0.28ex]}

\opn\dirlim{\underrightarrow{\lim}} \opn\inivlim{\underleftarrow{\lim}}
\let\union=\cup
\let\sect=\cap
\let\dirsum=\oplus
\let\tensor=\otimes
\let\iso=\cong
\let\Union=\bigcup
\let\Sect=\bigcap
\let\Dirsum=\bigoplus

\let\to=\rightarrow
\let\To=\longrightarrow
\def\Implies{\ifmmode\Longrightarrow \else
        \unskip${}\Longrightarrow{}$\ignorespaces\fi}
\def\implies{\ifmmode\Rightarrow \else
        \unskip${}\Rightarrow{}$\ignorespaces\fi}
\def\iff{\ifmmode\Longleftrightarrow \else
        \unskip${}\Longleftrightarrow{}$\ignorespaces\fi}

\let\:=\colon
\newtheorem{Theorem}{Theorem}[section]
\newtheorem{Lemma}[Theorem]{Lemma}
\newtheorem{Corollary}[Theorem]{Corollary}
\newtheorem{Proposition}[Theorem]{Proposition}
\newtheorem{Remark}[Theorem]{Remark}

\newtheorem{Example}[Theorem]{Example}
\newtheorem{Examples}[Theorem]{Examples}
\newtheorem{Definition}[Theorem]{Definition}

\let\epsilon\varepsilon
\let\phi=\varphi
\let\kappa=\varkappa
\def\qed{\ifhmode\textqed\fi
      \ifmmode\ifinner\quad\qedsymbol\else\dispqed\fi\fi}
\def\textqed{\unskip\nobreak\penalty50
       \hskip2em\hbox{}\nobreak\hfil\qedsymbol
       \parfillskip=0pt \finalhyphendemerits=0}
\def\dispqed{\rlap{\qquad\qedsymbol}}

\opn\dis{dis}
\def\pnt{{\raise0.5mm\hbox{\large\bf.}}}

\opn\Lex{Lex}

\begin{document}

\title{Finite filtrations of modules and shellable multicomplexes}

\author{J\"urgen Herzog and   Dorin Popescu}
\thanks{The second author was mainly supported by Marie Curie Intra-European
Fellowships MEIF-CT-2003-501046 and partially supported by CNCSIS and the Ceres programs 4-147/2004 and 4-131/2004 of the Romanian Ministery of Education and
Research}
\address{J\"urgen Herzog, Fachbereich Mathematik und
Informatik, Universit\"at Duisburg-Essen, Campus Essen, 45117 Essen, Germany} \email{juergen.herzog@uni-essen.de}

\address{Dorin Popescu, Institute of Mathematics "Simion Stoilow", University of Bucharest,
P.O.Box 1-764, Bucharest 014700, Romania} \email{dorin.popescu@imar.ro} \maketitle

\begin{abstract}
We introduce pretty clean modules, extending  the notion of clean modules  by Dress, and  show that pretty clean modules are sequentially
Cohen-Macaulay. We also extend a theorem of Dress on shellable simplicial complexes to multicomplexes.
\end{abstract}

\section*{Introduction}
Let $R$ be a Noetherian ring, and $M$ a finitely generated $R$-module.
 A basic fact in
commutative algebra (see \cite[Theorem 6.4]{M}) says that there exists a finite filtration $$\mathcal{F}\: 0=M_0\subset M_1\subset \cdots \subset M_{r-1}
\subset M_r=M$$ with cyclic quotients
 $M_i/M_{i-1}\iso R/P_i$ and   $P_i\in \Supp(M)$. We call any such filtration of $M$ a
 {\em prime filtration}. The set of prime ideals $\{P_1,\ldots, P_r\}$
 which define
 the cyclic quotients of $\mathcal{F}$ will be denoted by $\Supp(\mathcal{F})$. Another basic fact \cite[Theorem 6.5]{M} says that
\[
\Ass(M)\subset \Supp(\mathcal{F})\subset \Supp(M).
\]
Let $\Min(M)$ denote the set of minimal prime ideals. Dress \cite{D} calls a prime filtration ${\mathcal F}$ of $M$ {\em clean}, if
$\Supp(\mathcal{F})\subset\Min(M)$. The module  $M$ is called {\em clean}, if $M$ admits a clean filtration. It is clear that for a clean filtration $\mathcal
F$ of $M$ one has $$\Min(M)=\Ass(M)=\Supp({\mathcal F}).$$

Cleanness is the algebraic counterpart of shellability for simplicial complexes. Indeed, let $\Delta$ be a simplicial complex and $K$ a field. Dress \cite{D}
showed that $\Delta$ is (non-pure) shellable in the sense of Bj\"orner and Wachs \cite{BW}, if and only if the Stanley-Reisner ring $K[\Delta]$ is clean.

On the other hand Stanley \cite{St} showed that if $\Delta $ is  shellable, then $K[\Delta]$ is sequentially Cohen-Macaulay. In this paper we show more
generally that any clean module  over a Cohen-Macaulay ring which admits a canonical module is sequentially Cohen-Macaulay if  all factors in the clean
filtration are Cohen-Macaulay. In fact, we prove this result (Theorem \ref{sequentially}) for an even larger class of modules which we call pretty clean. These
modules are defined by the property that they have a prime filtration as above, and  such that for all $i<j$ for which $P_i \subset P_j$ it follows that
$P_i=P_j$.

We now describe the content of this paper in more detail. In Section 1 we recall the concept of dimension filtrations introduced by Schenzel \cite{Sc}, and
note (Proposition \ref{characterization}) that the dimension  filtration of a module is characterized  by the associated prime ideals of its factors.  In the
next section we discuss some basic properties of sequentially Cohen-Macaulay modules. Such modules were introduced by Schenzel \cite{Sc} and Stanley \cite{St}.
It was Schenzel who observed that a module is sequentially Cohen-Macaulay if and only the non-zero factors of the dimension filtration are Cohen-Macaulay.

The following section is devoted to introduce clean and pretty clean modules. We show that a pretty clean filtration $\mathcal F$ of a module $M$ satisfies
$\supp({\mathcal F})=\Ass(M)$, and we give an example of a module $M$ which admits a prime filtration ${\mathcal F}$ with $\supp({\mathcal F})=\Ass(M)$ but
which is not pretty clean. We also observe that that all pretty clean filtrations of a module have the same length.

In Section 4 we show (Theorem \ref{sequentially}) that under the mild assumptions, mentioned above, pretty clean modules are sequentially Cohen-Macaulay, and
we show in Corollary \ref{interesting} that under the same assumptions a module is pretty clean if and only if the factors in its dimension filtration are all
clean.

In Section 5 we give an interesting class of pretty clean rings, namely of rings whose defining ideal is of Borel type. This generalizes a result in \cite{HPV}
where it is shown that such rings are sequentially Cohen-Macaulay.

In the following section we consider graded and multigraded pretty clean rings and modules. Of particular interest is the case that $R=S/I$ where $S$ is a
polynomial ring and $I\subset S$ a monomial ideal. Using a result of Nagel and R\"omer \cite[Theorem 3.1]{NR} we show that in this case the length of each
multigraded pretty clean filtrations of $S/I$  is equals to the arithmetic degree of $S/I$.

In \cite{St1} Stanley conjectured that the depth of $S/I$ is a lower bound for the `size'  of the summands in any Stanley decomposition of $S/I$. We show in
Theorem \ref{stanley1} that Stanley's conjecture holds if $R$ is a multigraded pretty clean ring.

In Section 7 we show that for a given prime filtration $\mathcal{F}\: 0=M_0\subset M_1\subset \cdots \subset M_{r-1} \subset M_r=M$ of $M$ with factors
$M_i/M_{i-1}=R/P_i$ there exists irreducible submodules $P_j$-primary submodules $N_j$ of $M$ such that $M_i=\Sect_{j>i}^rN_j$ for $i=0,\ldots, r$. It turns
out, as demonstrated in the next and the following sections, that this presentation of the modules $M_i$ is the algebraic interpretation of shellability for
clean and pretty clean filtrations. This becomes obvious in the next section where we recall the theorem  of Dress and show that the shelling numbers of a
simplicial complex can be recovered from the graded clean filtration, see Proposition \ref{shelling numbers}.

In Section 9 we introduce multicomplexes. These are subsets  $\Gamma\subset \NN^n_\infty$ which are closed under limits of sequence $a_i\in \Gamma$  with
$a_i\leq a_{i+1}$ (componentwise), and have the property that whenever $a\in \Gamma$ and $b\leq a$ (componentwise), then $b\in \Gamma$. Here
$\NN_\infty=\NN\union \{\infty\}$. We show that if
 $\Gamma$ is a multicomplex and $a\in \Gamma$, then there exists a maximal element $m\in \Gamma$ with $a\leq m$. Here we need that $\Gamma$ is closed with respect
 to limits of non-decreasing sequences. Then we define the facets of $\Gamma$ to be those elements $a\in\Gamma$  with the property that if $a\leq m$ and $m$ is
 maximal in $\Gamma$, then the infinite part of $a$ coincides with the infinite part of $m$, which means that the $i$th component of $a$ is infinite if and only if the
 $i$th component of $m$ is infinite. We show that each multicomplex has only a finite number of facets.

 Multicomplexes in $\NN^n_\infty$ correspond to monomial ideals in $S=K[x_1,\ldots,x_n]$. The monomial ideal $I$ defined by a multicomplex  $\Gamma$ is the ideal  spanned by all
 monomials whose exponents belong to $\NN^n\setminus \Gamma$. Our definition of the facets of $\Gamma$ is partly justified by the fact, shown in
Lemma \ref{pairs}, that there  is a bijection between the set  of facets of $\Gamma$ and the standard pairs of $I$ as defined by Sturmfels, Trung and Vogel in
\cite{STV}. However the main justification of the definition is given by Proposition \ref{multiprimary} where we show that a pretty clean filtration of $S/I$
determines uniquely the facets of $\Gamma$. This result finally leads us to the definition of shellable multicomplexes. In Proposition \ref{extend} we show
that our definition of shellable multicomplexes extends the corresponding notion known for simplicial complexes. However the main result of the final section
is Theorem \ref{multi2} which asserts  that for a monomial ideal $I$ the ring $S/I$ is  multigraded pretty clean  if and only if the corresponding multicomplex
is shellable.

\section{The  dimension  filtration}

Let $M$ be an $R$-module of dimension $d$. In \cite{Sc} Schenzel introduced the {\em dimension filtration}
\[
{\mathcal F}\: 0\subset D_0(M)\subset D_1(M)\subset \cdots \subset D_{d-1}(M)\subset D_d(M)=M
\]
of $M$, which is defined by the property that $D_i(M)$ is the largest submodule of $M$ with $\dim D_i(M)\leq i$ for $i=0,\ldots,d$. It is convenient to set
$D_{-1}(M)=(0)$.

For all $i$ we set  $\Ass^i(M)=\{P\in\Ass(M)\: \dim R/P=i\}$. The following characterization of a dimension filtration will be useful for us:

\begin{Proposition}
\label{characterization} Let ${\mathcal F}\: 0\subset M_0\subset M_1\subset \cdots \subset M_{d-1}\subset M_d=M$ be a filtration of $M$. The following
conditions are equivalent:
\begin{enumerate}
\item[(a)] $\Ass(M_i/M_{i-1})=\Ass^{i}(M)$ for all $i$; \item[(b)] ${\mathcal F}$ is the dimension filtration of $M$.
\end{enumerate}
\end{Proposition}

\begin{proof}
That the dimension filtration satisfies condition (a) has been shown by Schenzel in \cite[Corollary 2.3 (c)]{Sc}.

For the converse we show that if ${\mathcal F}$ satisfies condition (a), then it is uniquely determined. Since the dimension filtration satisfies this
condition, it follows then that ${\mathcal F}$ must be the dimension filtration of $M$.

The integers $i$ for which $M_i=M_{i-1}$ are exactly those for which $\Ass^{i}(M)=\emptyset$, and hence this set is uniquely determined.

Thus it remains to show, if $M_i\neq M_{i+1}$, then $M_{i}$ is uniquely determined. To this end, consider the multiplicatively closed set
\[
S=R\setminus\Union_{P\in\Ass(M),\atop \dim R/P\geq i+1}P,
\]
and let $U$ be the kernel of the natural map $M\to M_S$. We claim that $M_i=U$. This will imply the uniqueness of the filtration.

We first notice that  $(M_j/M_{j-1})_S=0$ for $j\leq i$.  Indeed, if $(M_j/M_{j-1})_S\neq 0$, then  $PR_S\in\Ass_{R_S}(M_j/M_{j-1})_S$ for some   $P\in
\Ass_R(M_j/M_{j-1})$. By (a),  $\dim R/P\leq i$, and hence  $P\sect S\neq\emptyset$,  a contradiction. We conclude that $(M_i)_S=0$, and hence $M_i\subset U$.

Condition  (a) implies that
\[
\Ass(M/M_{i})\subset \Union_{P\in\Ass(M),\atop \dim R/P\geq  i+1}P.
\]
Therefore all elements of $S$ are non-zerodivisors on $M/M_{i}$, and hence the natural map $M/M_i\to (M/M_i)_S$ is injective. This implies that $U\subset M_i$.
\end{proof}

It follows from  condition  (a) of Proposition \ref{characterization} that if $M_i/M_{i-1}\neq 0$, then $M_i/M_{i-1}$ is equidimensional  of dimension $i$ and
has no embedded prime ideals.

The arguments in the proof of the previous proposition yield the following description of the dimension filtration.

\begin{Corollary}[Schenzel]
\label{description} Let $\Sect_{i=1}^nQ_i$ be a  primary decomposition of $(0)$ in $M$, where $Q_i$ is $P_i$-primary. Then
\[
D_i(M)=\Sect_{\dim R/P_j\geq i+1}Q_j,
\]
for $i=1,\ldots, \dim M$.
\end{Corollary}

\section{Sequentially Cohen-Macaulay modules}

Let $(R,\mm)$ be a local Noetherian ring, or a standard graded $K$-algebra with graded maximal ideal $\mm$. All modules considered will be finitely generated,
and graded if $R$ is graded.

The following definition is due to Stanley \cite[Section II, 3.9]{St}, and  Schenzel \cite{Sc}.

\begin{Definition}
\label{stanley} {\em Let $M$ be a finitely generated (graded) $R$-module. A finite filtration
\[
0=M_0\subset M_1\subset M_2\subset\ldots\subset M_r=M
\]
of $M$ by (graded) submodules of $M$ is called a {\em CM-filtration}, if   each quotient $M_i/M_{i-1}$ is Cohen-Macaulay (CM for short),   and  $$\dim
(M_1/M_0)<\dim(M_2/M_1)<\ldots <\dim(M_r/M_{r-1}).$$

The module $M$ is called {\em sequentially Cohen-Macaulay} if $M$ admits a CM-filtration. }

\end{Definition}

We recall a few basic facts whose proof in the graded case can be found in \cite{HS}, but which  are proved word by word in the same way  in the local case.

\begin{Proposition}
\label{exti} Let $R$ be Cohen-Macaulay of dimension $n$ with canonical module $\omega_R$. Suppose that $M$ is sequentially CM with a CM-filtration as in {\em
\ref{stanley}}, and assume further that $d_i=\dim M_i/M_{i-1}$ for $i=1,\ldots, r$. Then
\begin{enumerate}
\item[(a)]  $\Ext^{n-d_i}_R(M,\omega_R)\iso \Ext^{n-d_i}_R(M_i/M_{i-1},\omega_R)$; \item[(b)] $\Ext^{n-d_i}_R(M,\omega_R)$ is CM of dimension $d_i$ for
$i=1,\ldots,r$; \item[(c)] $\Ext^j_R(M,\omega_R)=0$ if $j\not\in\{n-d_1,\ldots,n-d_r\}$; \item[(d)] $\Ext^{n-d_i}_R(\Ext^{n-d_i}_R(M,\omega_R),\omega_R)\iso
M_i/M_{i-1}$ for $i=1,\ldots,r$.
\end{enumerate}
\end{Proposition}

\begin{Corollary}
\label{Ass} With the assumptions and notation introduced in Proposition \ref{exti} we have
\[
\Ass(\Ext^{n-d_i}_R(M,\omega_R))=\Ass(M_i/M_{i-1}).
\]
\end{Corollary}

\begin{proof}
Let ${\bf x}=x_1,\ldots, x_{n-d_i}$ be a maximal regular sequence in $\Ann(M_i/M_{i-1})$, and set $S=R/({\bf x})$. Then (a) implies that
$\Ext^{n-d_i}_R(M,\omega_R)\iso \Hom_S(M_i/M_{i-1},\omega_S)$, and that $M_i/M_{i-1}$ may be viewed a maximal CM module over $S$.  It follows that
\[
\Ass(\Ext^{n-d_i}_R(M,\omega_R))=\Supp(M_i/M_{i-1})\sect\Ass(\omega_S)=\Supp(M_i/M_{i-1})\sect \Min(S).
\]
Since $M_i/M_{i-1}$ is a maximal CM module over $S$, we have
\[
\Ass(M_i/M_{i-1})=\Min(M_i/M_{i-1})=\Supp(M_i/M_{i-1})\sect \Min(S).
\]
This proves the assertion.
\end{proof}

 It follows from Proposition \ref{exti}  that if $M$ is sequentially CM, then the  non-zero modules among the $\Ext_R^{n-i}(M,\omega_R)$ are CM of dimension $i$. Peskine noticed that this property characterizes sequentially CM modules. Indeed one has

\begin{Theorem}
\label{peskine} The following two conditions are equivalent:
\begin{enumerate}
\item[(a)] $M$ is sequentially CM;
\item[(b)] for all $i$, the modules $\Ext_R^{n-i}(M,\omega_R)$ are either $0$ or CM of 
dimension $i$.
\end{enumerate}
\end{Theorem}

We conclude this section with a result due to  Schenzel \cite[Corollary 2.3]{Sc}.

\begin{Proposition}
\label{faridi} Let $M$ be sequentially CM with a CM-filtration as above. Then $\Ass(M_i/M_{i-1})=\Ass^{d_i}(M)$ for all $i$. In particular,
$\Ass(M)=\Union_i\Ass(M_i/M_{i-1})$.
\end{Proposition}

\begin{proof}
Since $M_i/M_{i-1}$ is CM of dimension $d_i$, it follows that $\dim R/P=d_i$ for all $P\in \Ass(M_i/M_{i-1})$. Therefore it suffices to show that
$\Ass(M)=\Union_i\Ass(M_i/M_{i-1})$. Using the fact that for an exact  sequence $0\to U\to V\to W\to 0$ of $R$-modules one has that $\Ass(V)\subset
\Ass(U)\union \Ass(W)$, one easily concludes by induction on the length $r$ of the filtration, that $\Ass(M)\subset  \Union_i\Ass(M_i/M_{i-1})$.

Conversely, let $P\in \Ass(M_i/M_{i-1})$. Then $(M_j/M_{j-1})_P=0$ for all $j<i$, since $\dim M_j/M_{j-1}<\dim R/P$. This implies that $(M_{i-1})_P=0$, so that
$(M_i/M_{i-1})_P=(M_i)_P$. Thus $PR_P\in \Ass_{R_P}(M_i)_P$, and hence $P\in \Ass(M_i)\subset \Ass(M)$.
\end{proof}

Combining Proposition \ref{faridi} with Corollary \ref{Ass} we obtain

\begin{Corollary}
\label{need} Let $M$ be sequentially CM, then $\Ass(\Ext_R^{n-i}(M,\omega_R))=\Ass^i(M)$ for all $i$.
\end{Corollary}

The following characterization of sequentially Cohen-Macaulay modules, due to Schenzel \cite[Proposition 4.3]{Sc},  is a consequence of Proposition
\ref{characterization} and  Proposition \ref{faridi}.

\begin{Corollary}
\label{comparison} A module  $M$ is sequentially CM, if and only if the factors in the dimension filtration of $M$ are either 0 or CM.
\end{Corollary}

\section{Clean and pretty clean modules}

Let $R$ be a Noetherian ring, and $M$ a finitely generated $R$-module. Recall from the introduction that according to Dress \cite{D} a prime filtration
$\mathcal F$ of $M$ is called {\em clean} if $\Supp({\mathcal F})=\Min(M)$, and that $M$ itself is called {\em clean} if $M$ admits a clean filtration.

\begin{Lemma}
\label{character} Let $\mathcal F$ be a prime filtration of $M$. The following conditions are equivalent:
\begin{enumerate}
\item[(a)] $\mathcal F$ is a clean filtration of $M$; \item[(b)] For all  $P, Q\in \Supp({\mathcal F})$ with  $P\subset Q$ one has $P=Q$.
\end{enumerate}
\end{Lemma}

\begin{proof}
(a)\implies (b) is trivial. Conversely suppose ${\mathcal F}: 0=M_0\subset M_1\subset\cdots\subset M_r=M$ with $M_i/M_{i-1}=R/P_i$ and let  $P\in
\Supp({\mathcal F})$. Since there are no non-trivial inclusions between the prime ideals in $\Supp({\mathcal F})$ it follows that $M_P$ has a filtration
$(0)=(M_0)_P\subset (M_1)_P\subset \cdots \subset (M_r)_P=M_P$ such that
$$(M_i)_P/(M_{i-1})_P = \left\{ \begin{array}{lll} R_P/PR_P, & \text{if} & P=P_i,\\
0, & \mbox{if} & P\neq P_i. \end{array} \right.
$$
Hence we see that $\Ass_{R_P}(M_P)=\{PR_P\}$, and so $P\in \Ass(M)$. It follows that $\Supp({\mathcal F})=\Ass(M)$. Applying again  assumption (b), we conclude
that $\Ass(M)=\Min(M)$.
\end{proof}

\begin{Corollary}
\label{induced} Let $0=M_0\subset M_1\subset \ldots \subset M_{r-1}\subset M_r=M$ be a clean filtration of $M$. Then for all $i=0,\ldots,r$
\[
0=M_i/M_i\subset M_{i+1}/M_i\subset\ldots \subset  M_{r-1}/M_i\subset M_r/M_i,
\]
and
\[
0=M_0\subset M_1\subset \ldots\subset M_{i-1}\subset M_i
\]
are clean filtrations. In particular, $M_i$ and $M/M_i$ are clean.
\end{Corollary}

A weakening of  condition (b) of Lemma \ref{character} leads to

\begin{Definition}
\label{pretty} {\em  A prime filtration ${\mathcal F}\: 0=M_0\subset M_1\subset \ldots \subset M_{r-1}\subset M_r=M$ of $M$ with $M_i/M_{i-1}=R/P_i$ is called
{\em pretty clean}, if for all $i<j$ for which  $P_i\subset P_j$ it follows that $P_i=P_j$.

In other words, a proper inclusion $P_i\subset P_j$ is only possible if $i>j$. The module $M$ is called {\em pretty clean}, if it has a pretty clean
filtration. A ring is called pretty clean if it is a pretty clean module, viewed as a module over itself.

}
\end{Definition}

\begin{Remark}
{\em Let ${\mathcal F}\: 0=M_0\subset M_1\subset \ldots \subset M_{r-1}\subset M_r=M$ be a pretty clean filtration of $M$. It follows immediately from the
definition that for all $i$ the filtrations
\[
0=M_i/M_i\subset M_{i+1}/M_i\subset\ldots \subset  M_{r-1}/M_i\subset M_r/M_i,
\]
and
\[
0=M_0\subset M_1\subset \ldots\subset M_{i-1}\subset M_i
\]
are pretty clean. }
\end{Remark}

\begin{Lemma}
\label{important} Let  ${\mathcal F}\: 0=M_0\subset M_1\subset \ldots \subset M_{r-1}\subset M_r=M$ be a pretty clean filtration of $M$. Then $P_i\in
\Ass(M_i)$ for all $i$.
\end{Lemma}

\begin{proof}
We use the same  argument as in the proof of Lemma \ref{character}: set $P=P_i$. Then $0=(M_0)_P\subset (M_1)_P\subset \ldots \subset (M_{i-1})_P\subset
(M_{i})_P$ is a finite  filtration of the $R_P$-module $(M_i)_P$. Let $j\leq i$. Since ${\mathcal F}$ is pretty clean we get
$$(M_j)_P/(M_{j-1})_P = \left\{ \begin{array}{lll} R_P/PR_P, & \text{if} & P_j=P,\\
0, & \mbox{if} & P_j\neq P. \end{array} \right.
$$
This implies that $PR_P\in \Ass_{R_p}((M_i)_P)$. Therefore $P\in \Ass(M_i)$.
\end{proof}

\begin{Corollary}
\label{nice} Let ${\mathcal F}$ be a pretty clean filtration of $M$. Then $\Supp({\mathcal F})=\Ass(M)$.
\end{Corollary}

\begin{proof}
For all $i$ we have $P_i\in \Ass(M_i)\subset \Ass(M)$. Therefore $\Supp({\mathcal F})\subset \Ass(M)$. The other inclusion holds for any prime filtration.
\end{proof}

\begin{Corollary}
\label{equal} Let $M$ be a pretty clean module. The following conditions are equivalent:
\begin{enumerate}
\item[(a)] $M$ is clean; \item[(b)] $\Ass(M)=\Min(M)$.
\end{enumerate}
\end{Corollary}

\begin{Examples}{\em   Let $S=K[x,y]$ be the polynomial ring over the field $K$, $I\subset S$ the ideal $I=(x^2,xy)$ and $R=S/I$. Then $R$ is pretty clean but not clean.
Indeed, $0\subset (x)\subset R$ is a pretty clean filtration of $R$ with $(x)=R/(x,y)$, so that $P_1=(x,y)$ and $P_2=(x)$. $R$ is not clean since $\Ass(R)\neq
\Min(R)$.

Note $R$ has a different prime filtration, namely, ${\mathcal G}: 0\subset (y)\subset (x,y)\subset R$ with factors $(y)=R/(x)$ and $(x,y)/(y)=R/(x,y)$. Hence
this filtration is not pretty clean, even though $\Supp({\mathcal G})=\Ass(M)$. On the other hand, in the next section we give an example of a module which
admits a prime filtration ${\mathcal F}$ with $\supp({\mathcal F})=\Ass(M)$, but which is not pretty clean.}
\end{Examples}

We conclude this section by showing that all pretty clean filtrations have the same length. For $\pp\in\Spec(R)$ the number
\[
\mult_M(\pp)=\length (H^0_\pp(M_\pp)),
\]
is called the {\em length multiplicity} of $\pp$ with respect to $M$. Obviously, one has $\mult_M(\pp)>0$, if and only if $\pp\in \Ass(M)$. Localizing a pretty
clean filtration of $M$ we immediately get

\begin{Lemma}
\label{allthesame} Let $M$ be a pretty clean module. Then all pretty clean filtrations of $M$ have the same length, namely their common  length equals
$\sum_{\pp\in\Ass(M)}\mult_M(\pp)$.
\end{Lemma}

Assume now that $(R,\mm)$ is local.  Recall that the {\em arithmetic degree} of $M$ is defined to be $\sum_\pp\mult_M(\pp)\deg(R/\pp)$ where $\deg(R/\pp)$ is the multiplicity of the associated graded ring of $R/\pp$. The preceding lemma shows that the length of a pretty clean
filtration is bounded above by the arithmetic degree of the module, and equals the arithmetic degree if and only if $\deg R/\pp=1$ for all $\pp\in\Ass(M)$.

\section{Pretty clean modules are sequentially Cohen-Macaulay}

The purpose of this section is to show

\begin{Theorem}
\label{sequentially} Let $R$ be a local or standard graded CM ring admitting a canonical module $\omega_R$, and let $M$ be an $R$-module with pretty clean
filtration ${\mathcal F}$ such that $R/P$ is Cohen-Macaulay for all $P\in \Supp({\mathcal F})$. Furthermore suppose that $M$ is graded if $R$ is graded. Then
$M$ is sequentially Cohen-Macaulay. Moreover, if $\dim R/P=\dim M$ for all $P\in\Supp({\mathcal F})$, then $M$ is clean and Cohen-Macaulay.
\end{Theorem}

\begin{proof}
Let $n=\dim R$. We may assume that $R$ is local. In the graded case the arguments are the same.

For all $i$ we show: the module  $\Ext_R^{n-i}(M,\omega_R)$ is either $0$ or Cohen-Macaulay of dimension $i$. We show this by induction on the length $r$ of
the pretty clean filtration ${\mathcal F}\: 0=M_0\subset M_1\subset\cdots \subset M_{r-1}\subset M_r=M$ of $M$. Since, as we already noticed, the module
$U=M_{r-1}$  is pretty clean with a pretty clean filtration of length $r-1$,  we may assume by induction that $U$ is sequentially Cohen-Macaulay.

Let $M/U=R/P$. By hypothesis, $R/P$ is Cohen-Macaulay, say of dimension $d$. The short exact sequence
\[
0\To U\To M\To R/P\To 0
\]
gives rise to the long exact sequence
\[
\cdots \Ext_R^{n-i-1}(U,\omega_R)\to \Ext_R^{n-i}(R/P,\omega_R)\to \Ext_R^{n-i}(M,\omega_R)\to\Ext_R^{n-i}(U,\omega_R)\to\cdots,
\]
Since
\[
\Ext_R^{n-i}(R/P,\omega_R)= \left\{ \begin{array}{lll} \omega_{R/P}, & \text{if} & i=d\,\\
0, & \mbox{if} & i\neq d, \end{array} \right.
\]
it follows that $\Ext_R^{n-i}(M,\omega_R)\iso \Ext_R^{n-i}(U,\omega_R)$ for all $i\neq d,d+1$. Thus for such $i$ we have $\Ext_R^{n-i}(M,\omega_R)$ is
Cohen-Macaulay of dimension $i$ if not the zero module.

Moreover we have the exact sequence
\begin{eqnarray*} 0&\to &\Ext_R^{n-d-1}(M,\omega_R)\to \Ext_R^{n-d-1}(U,\omega_R)\to \Ext_R^{n-d}(R/P,\omega_R)\\
 &\to& \Ext_R^{n-d}(M,\omega_R)\to
\Ext_R^{n-d}(U,\omega_R)\to 0.
\end{eqnarray*}

Suppose the map $\Ext_R^{n-d-1}(U,\omega_R)\to \Ext_R^{n-d}(R/P,\omega_R)\iso\omega_{R/P}$ is not the zero map. Then its image $C\subset \omega_{R/P}$ is not
zero. Since $R/P$ is domain, $\omega_{R/P}$ may be identified with an ideal in $R/P$, see \cite[Proposition 3.3.18]{BH}. Hence also $C$ may be identified with
an ideal in $R/P$. Again using that $R/P$ is a domain, we conclude that $CR_P\neq 0$. It follows that $\Ext_R^{n-d-1}(U,\omega_R)_P\neq 0$, and so the set
\[
\Ass_{R_P}(\Ext_R^{n-d-1}(U,\omega_R)_P)=\{QR_P\: Q\in \Ass_R(\Ext_R^{n-d-1}(U,\omega_R)),\; Q\subset P\}
\]
is not empty. Thus there exists  $Q\in \Ass_R(\Ext_R^{n-d-1}(U,\omega_R))$ with $Q\subset P$. By Corollary \ref{need} we know that
$\Ass_R(\Ext_R^{n-d-1}(U,\omega_R))\subset \Ass_R^{d+1} (U)$. Therefore, since $\dim R/P=d$, the inclusion $Q\subset P$ must be proper. But this contradicts
the fact that ${\mathcal F}$ is a pretty clean filtration of $M$.

It follows now that
\[
\Ext_R^{n-d-1}(M,\omega_R)\iso  \Ext_R^{n-d-1}(U,\omega_R),
\]
and that   the  sequence
\begin{eqnarray}
\label{exact}
 0\To\omega_{R/P}\to \Ext_R^{n-d}(M,\omega_R)\To \Ext_R^{n-d}(U,\omega_R)\To 0
\end{eqnarray}
is exact. Using the induction hypothesis we conclude  that $\Ext_R^{n-d-1}(M,\omega_R)$ is either Cohen-Macaulay  of dimension $d+1$ or the zero module, and
that $\Ext_R^{n-d}(M,\omega_R)$ is Cohen-Macaulay of dimension $d$.

If $\dim R/P=\dim M$ for all $P\in \Supp({\mathcal F})$, then the pretty clean filtration ${\mathcal F}$ is necessarily clean, and $M$ is unmixed. Since any
unmixed sequentially Cohen-Macaulay module is Cohen-Macaulay, all assertions are proved.
\end{proof}

As a consequence of the previous theorem we get

\begin{Corollary}
\label{interesting} Let $M$ be an $R$-module. If the  non-zero factors of the dimension filtration of $M$ are clean, then $M$ is pretty clean.

Conversely assume that $R$ is a local or standard graded  CM ring with canonical module $\omega_R$, and that  $M$ admits a pretty clean filtration  ${\mathcal
F}$ such that $R/P$ is CM for all $P\in \Supp({\mathcal F})$. Furthermore assume that $M$ is graded if $R$ is graded. Then the  non-zero factors of the
dimension filtration of $M$ are clean.
\end{Corollary}

\begin{proof} Suppose  all factors $D_i(M)/D_{i-1}(M)$ in the dimension filtration of $M$ are clean. Then it is obvious that the dimension filtration can be refined
to yield a pretty clean filtration of $M$.

We prove the second statement  of the corollary by induction on the length $r$ of the filtration $\mathcal F$. The claim is obvious if $r=1$. Now let $r>1$,
and set $U=M_{r-1}$. We obtain the  exact sequence $0\to U\to M\to R/P\to 0$ with $P\in \Spec(R)$.  Let $d=\dim R/P$. Then, as we have seen in the proof of
Theorem \ref{sequentially}, one has $ \Ext_R^{n-i}(M,\omega_R)\iso \Ext_R^{n-i}(U,\omega_R)$ for all $i\neq d$, as well as the exact sequence
\[
0\To\omega_{R/P}\to \Ext_R^{n-d}(M,\omega_R)\To \Ext_R^{n-d}(U,\omega_R)\To 0.
\]
Since $M$ is sequentially CM by the previous theorem, these isomorphisms together with Proposition \ref{exti}(d) and Corollary \ref{comparison} imply that
$$D_i(M)/D_{i-1}(M)\iso D_i(U)/D_{i-1}(U)$$  for $i\neq d$. Hence, since the factors $D_i(U)/D_{i-1}(U)$ are clean by induction hypothesis, the same is true for
the factors $D_i(M)/D_{i-1}(M)$ with $i\neq d$.

Applying the functor $\Ext_R^{n-d}(-,\omega_R)$ to the above exact sequence and using Proposition \ref{exti}(d) again  we obtain the exact sequence
\[
0\To D_d(U)/D_{d-1}(U)\To D_d(M)/D_{d-1}(M)\To R/P\to 0.
\]
Since all modules in this exact sequence are of dimension $d$, and since $D_d(U)/D_{d-1}(U)$ is clean, it follows that $D_d(M)/D_{d-1}(M)$ is clean as well.
\end{proof}

\begin{Corollary}
\label{conclusion}
 Let $S=K[x_1,\ldots,x_n]$ be a the polynomial ring and $I\subset S$  a monomial ideal. Then the following conditions are equivalent:
 \begin{enumerate}
 \item[(a)] $S/I$ is pretty clean;
 \item[(b)] $S/I$ is sequentially CM, and the non-zero factors in the dimension filtration of $S/I$ are clean;
\item[(c)] the non-zero factors in the dimension filtration of $S/I$ are clean.

 \end{enumerate}
 \end{Corollary}

 \begin{proof}
 (a)\implies (b):  Since the associated prime ideals of $S/I$ are all generated by subsets of $\{x_1,\ldots,x_n\}$, all hypotheses Theorem
 \ref{sequentially} and Corollary \ref{interesting} are satisfied, so that the assertions follow.

(b)\implies (c) is trivial.

 (c)\implies (a): The refinement of the dimension filtration by the clean filtrations of the non-zero factors gives us the desired pretty clean filtration of
$S/I$.
\end{proof}

\begin{Example}{\em Let $S=K[x,z,u,v]$, and consider the ideals $L=(u,v,z)$, $Q_1=(x,z^2)$,  $Q_2=(x,v^2,z^3)$ and $I=L\sect Q_1\sect Q_2$. We claim that
the module  $M=L/I$ is not pretty clean, but that $M$ has a prime filtration $\mathcal F$ with $\Supp({\mathcal F})=\Ass(M)$.

Note that $(L\sect Q_1)\sect (L\sect Q_2)$ modulo $I$ is an irredundant primary decomposition of $(0)$ in $M$. Hence since $\emptyset \neq \Ass(L/L\sect
Q_i)\subset \Ass(S/Q_i)=\{P_i\}$ with $P_1=(x,z)$ and $P_2=(x,v,z)$ we see that $\Ass(M)=\{P_1,P_2\}$.

It follows from Corollary \ref{description} that $D_1(M)=(L\sect Q_1)/I$ and that $D_2(M)=M$. We show that $D_2(M)/D_1(M)=L/L\sect Q_1$ is not clean. Indeed,
suppose  $L/L\sect Q_1$ is clean. Then, since $\Ass(L/L\sect Q_1)=\{P_1\}$,  this module has a filtration with all factors isomorphic to $S/P_1$, and the
number of these factors equals the length of the $S_{P_1}$-module $(L/L\sect Q_1)_{P_1}=S_{P_1}/Q_1S_{P_1}$. This length is obviously $2$. On the other hand,
since $L/L\sect Q_1$ is generated by 3 elements, it cannot have a filtration with two factors, both of them being cyclic.

Knowing now that $D_2(M)/D_1(M)$ is not clean,  we conclude from Corollary \ref{interesting}  that $M$ is not pretty clean.

Finally we construct a prime filtration ${\mathcal F}$ of $M$ with $\Supp({\mathcal F})=\Ass(M)$. The filtration $\mathcal F$ will be the following  refinement
of the dimension filtration. Denote by $\bar{a}$ the residue class of an element $a\in L$ in $L/L\sect Q_1=D_2(M)/D_1(M)$. Then $(0)\subset
(\bar{z})\subset(\bar{z},\bar{v})\subset (\bar{z},\bar{v},\bar{u})=D_2(M)/D_1(M)$ is a filtration of $D_2(M)/D_1(M)$ with $(\bar{z})=S/P_1$,
$(\bar{z},\bar{v})/(\bar{z})=S/P_1$ and $(\bar{z},\bar{v},\bar{u})/(\bar{z},\bar{v})=S/P_2$. Furthermore, denote the residue class of an element $a\in S$ in
$S/I$ by $\tilde{a}$. Then $D_1(M)=L\sect Q_1/I$ is generated by $\tilde{z}^2$, and so $D_1(M)\iso S/(x,z,v^2)$. It is clear that this filtration can be further refined so
that all factors are isomorphic to $S/P_2$.
 }
\end{Example}

\section{Classes of pretty clean rings}

Let $S=K[x_1,\ldots, x_n]$ be the polynomial ring over a field $K$. In this section we present a class of monomial ideals for which $S/I$ is pretty clean.
Quite generally we have

\begin{Proposition}
\label{total} Let $I\subset S$ be a monomial ideal, and suppose that $\Ass(S/I)$ is totally ordered by inclusion. Then $S/I$ is pretty clean.
\end{Proposition}

\begin{proof}  Let $\Ass(S/I)=\{P_1,\ldots,P_r\}$ and suppose that $P_1\supset P_2\supset \cdots \supset P_r$, and set $d_i=\dim S/P_i$ for $i=1,\ldots, r$.
The ideal $I$  can be written as an intersection  $I=\Sect_{i=1}^rQ_i$  where each $Q_i$ is a $P_i$-primary monomial ideal. There exist subsets $J_i\subset
[n]$ such that  $P_i$ is generated by $x_j$ with $j\in J_i$. It follows from our assumption that $J_1\supset J_2\supset \cdots \supset J_r$.

Set $U_i=\Sect_{j>i}Q_j$. Then according to Corollary \ref{description} we have $U_i/I=D_{d_i}(S/I)$.  By Corollary \ref{conclusion}(c) it suffices to show
that $U_i/Q_i\sect U_i$ is clean for all $i$. We have $\emptyset \neq \Ass(U_i/Q_i\sect U_i)\subset \Ass(S/Q_i)=\{P_i\}$.  Let $S'$ be the polynomial ring over
$K$ in the variables $x_j$ with $j\in J_i$,  and set $P_i'=P_i\sect S'$. Then $P_i'$ is the graded maximal ideal of $S'$ and $P_i=P_i'S$. Similarly, since
$J_k\subset J_i$ for $k\geq i$, we have $Q_k=Q_k'S$ and $U_k=U_k'S$ where $Q_k'=Q_k\sect S'$ and $U_k'=U_k\sect S'$. The $S'$-module $U'_i/Q'_i\sect U'_i$ is a
clean since it is of finite length. By base change, $U_i/Q_i\sect U_i\iso (U'_i/Q'_i\sect U'_i)\tensor_{S'}S$ is a clean $S$-module.
\end{proof}

 In
Gr\"obner basis theory, Borel fixed ideals play an important role as they are just the generic initial ideals of graded ideals in a polynomial ring. By a
theorem of Bayer and Stillman (see \cite[Proposition 15.24]{Ei}) a Borel fixed ideal $I\subset S=K[x_1,\ldots,x_n]$ is a monomial ideal with the property that
\begin{eqnarray}
\label{borel}
 I:x_j^\infty=I:(x_1,\ldots,x_j)^\infty
\end{eqnarray}
for all $j=1,\ldots,n$. In \cite{HPV}, any monomial ideal satisfying condition (\ref{borel}) is called an ideal of {\em Borel type}, and it is shown   that
$S/I$ is sequentially Cohen-Macaulay if $I$ is of Borel type.

Here we show the following slightly stronger

\begin{Proposition}
\label{cleanborel} Let $I\subset S$ be an ideal of Borel type. Then $S/I$ is pretty clean.
\end{Proposition}

\begin{proof}
Let $P\in \Ass(S/I)$, and let $j$ be the largest integer such that $x_j\in P$. There exists  a monomial $u\in S$ such that $(I,u)/I\iso S/P$. Since $x_ju\in I$ it
follows that $u\in I:x_j^{\infty}$, and hence $u\in I:(x_1,\ldots,x_j)^{\infty}$. Therefore $u(x_1,\ldots, x_j)^k\subset I$ for some integer $k>0$,and hence
$(x_1,\ldots, x_j)^k\subset P$. Since $P$ is  prime ideal we conclude that $(x_1,\ldots, x_j)\subset P$. By the definition of $j$, it follows then that
$P=(x_1,\ldots, x_j)$.

Thus the associated prime ideal of $S/I$ are totally ordered and the assertion follows from Proposition \ref{total}.
\end{proof}

\section{Graded pretty clean modules}

Let $K$ be a  field and $R$  a standard graded $K$-algebra,  and let $M$ be a graded $R$-module. A prime filtration  of $M$
\[
{\mathcal F}\: (0)=M_0\subset M_1\subset M_{r-1}\subset M_r=M.
\]
is called {\em graded}, if all $M_i$ of $M$ are graded submodules of $M$, and if there  are graded isomorphisms  $M_i/M_{i-1}\iso R/P_i(-a_i)$ with some
$a_i\in\ZZ$ and some graded prime ideals $P_i$.

The module $M$ is called a {\em graded (pretty) clean module}, if it admits a (pretty) clean filtration which is a graded prime filtration.

Similarly we define multigraded  filtrations and multigraded (pretty) clean modules.

We denote by $(N)_i$ the $i$th graded component of a graded $R$-module $N$, and by
\[
\Hilb(N)=\sum_i\dim_K(N)_it^i\in\ZZ[t,t^{-1}]
\]
its Hilbert-series.

By the additivity of the Hilbert-series, one obtains for a  module with a graded prime filtration as above the Hilbert-series
\[
\Hilb(M)=\sum_{i=1}^r\Hilb(R/P_i)t^{a_i}.
\]

We now consider a more specific case

\begin{Proposition}
\label{monomial} Let $S=K[x_1,\ldots, x_n]$ be the polynomial ring, and $I\subset S$ a monomial ideal. Assume that $S/I$ is a graded pretty clean ring whose
graded pretty clean
 filtration has the  factors $M_j/M_{j-1}\iso S/P_j(-a_j)$  for $j=1,\ldots,
 r$, $a_j\in \NN$
 and  $P_j\in\Ass(S/I)$. For all $k$ and $i$ set
\[
h_{ki}=|\{j\: a_j=k,\; \dim S/P_j=i\}|.
\]
Then
\[
\Hilb(S/I)=\sum_iH_i(t) \quad \text{with}\quad H_i(t)=\frac{Q_i(t)}{(1-t)^i} \quad\text{where} \quad Q_i(t)=\sum_kh_{ki}t^k.
\]
\end{Proposition}

\begin{proof}
We have
\begin{eqnarray*}
\Hilb(S/I)&=&\sum_i \sum_{j\atop \dim S/P_j=i}\Hilb(S/P_j)t^{a_j}\\
&=& \sum_i(\sum_{j\atop \dim S/P_j=i}t^{a_j})/(1-t)^i.
\end{eqnarray*}
The last equality holds, since all associated prime ideals of $S/I$ are generated by subsets of the variables. Finally the desired formula follows, if we
combine in the sum $\sum_{j,\; \dim S/P_j=i}t^{a_j}$ all powers of $t$ with the same exponent.
\end{proof}

The attentive reader will notice the similarity of formula \ref{monomial} with the formula of McMullen and Walkup for shellable simplicial complexes, see
\cite[Corollary 5.1.14]{BH}. The precise relationship will become apparent in Section 8 where the numbers $a_j$ are interpreted as shelling numbers.

We now derive  similar formulas for the modules $\Ext^i_S(M,\omega_S)$ when $M$ is a graded  pretty clean module. Suppose $\dim S=n$. Using the graded version
of the exact sequence (\ref{exact}) in the proof of Theorem \ref{sequentially}, and induction on the length of the pretty clean filtration it follows easily
that
\[
\Hilb(\Ext_S^i(M,\omega_S))=\sum_{j\atop \dim S/P_j=n-i}\Hilb(\omega_{S/P_j})t^{-a_j} \quad\text{for}\quad i=0,\ldots, \dim M.
\]

In particular we have

\begin{Proposition}
\label{hilbert} With the assumptions and notation of {\em \ref{monomial}}, one has
\[
\Hilb(\Ext_S^i(S/I,S(-n)))=(\sum_k h_{k,n-i}t^{n-i-k})/(1-t)^{n-i}= (-1)^{n-i}H_{n-i}(t^{-1}).
\]
\end{Proposition}

\begin{proof}
The first equality follows from the fact that $\omega_{S/P_j}=S/P_j(-(n-i))$ if $\dim S/P_j=n-i$, so that $\Hilb(\omega_{S/P_j})=t^{n-i}/(1-t)^{n-i}$. To
obtain the second equality, we divide numerator and denominator of $(\sum_k h_{k,n-i}t^{n-i-k})/(1-t)^{n-i}$ by $t^{n-i}$ and get
\[
(\sum_k h_{k,n-i}t^{n-i-k})/(1-t)^{n-i}=(\sum_k h_{k,n-i}t^{-k})/(t^{-1}-1)^{n-i}=(-1)^{n-i}H_{n-i}(t^{-1}).
\]
\end{proof}

Let $S=K[x_1,\ldots, x_n]$ and $M$ be a graded $S$-module. We  set
$$b_j=\min\{k\: \Ext_S^j(M, S(-n))_k\neq 0\}.$$ Then the regularity
of $M$ is given by
$$\reg(M)=\max\{n-j-b_j\: j=0,\ldots, \},$$ cf.\ \cite[Section 20.5]{Ei}.

\begin{Corollary}
\label{regularity} Let $S=K[x_1,\ldots, x_n]$ be the polynomial ring, and $I\subset S$ a monomial ideal. Assume that $S/I$ is a graded pretty clean ring with
filtration as in {\em \ref{monomial}}. Then
\begin{enumerate}
\item[(a)] $\reg(S/I)=\max\{k\: h_{ki}\neq 0 \ \ \mbox{for some i}\}=\max\{a_j\: j=0,\ldots r\}$; \item[(b)] $\Hilb(D_i(S/I)/D_{i-1}(S/I))=H_i(t)$ for all $i$.
\end{enumerate}
\end{Corollary}

\begin{proof}
(a) the first equality follows immediately from Proposition \ref{hilbert} and the definition of $\reg(S/I)$. The second equality results from  the definition
of the numbers $h_{ki}$.

(b) By Proposition \ref{exti} we have
\[
D_i(S/I)/D_{i-1}(S/I)\iso\Ext_S^{n-i}(\Ext_S^{n-i}(S/I,\omega_S),\omega_S).
\]
Thus the assertion follows from Proposition \ref{hilbert} and \cite[Theorem 4.4.5(a)]{BH}.
\end{proof}

We denote by $e(M)$ the multiplicity of  a graded module.

\begin{Corollary}
\label{independent} Let $i$ and $k$ be integers. Then the number of factors $S/P(-k)$  in a graded pretty clean filtration of $S/I$ satisfying  $\dim S/P=i$ is
independent of the chosen filtration. In particular, all graded pretty clean filtrations of $S/I$ have the same length, namely
$\sum_{i=0}^ne(\Ext_S^i(S/I,S))$, and this number equals the arithmetic degree  of $S/I$.
\end{Corollary}

\begin{proof}
The  number in question equals $h_{ki}$, the $k$-th coefficient of the $h$-vector of $D_i(S/I)/D_{i-1}(S/I)$. Hence this number only depends on $S/I$.
Moreover, it follows that the length of a graded pretty clean filtration of $S/I$ equals $\sum_{i=0}^nQ_i(1)$. As a consequence of Proposition \ref{hilbert}
and  \cite[Proposition 4.1.9]{BH} we have that $Q_i(1)=e(\Ext_S^i(S/I,S))$. In \cite[Theorem 3.11]{NR}  Nagel and R\"omer have shown that the arithmetic degree
of a sequentially Cohen-Macaulay  $R$-module $M$ equals the number $\sum_{i=0}^ne(\Ext_R^i(M,\omega_R))$. Since  by Theorem \ref{sequentially}, $S/I$ is
sequentially CM, all assertions follow.
\end{proof}

We would like to remark that the fact that the length of all pretty clean filtrations of $S/I$ have length equal to the arithmetic degree of $S/I$ also follows
from Lemma \ref{allthesame}.

\medskip
\noindent
 Suppose that $I\subset S$ is a monomial ideal, and that $\mathcal F$ is a multigraded prime filtration of $S/I$ with factors $(S/P_i)(-a_i)$, $i=1,\ldots,r$,  where
$a_i\in \NN^n$. Then this filtration decomposes $S/I$ as a multigraded $K$-vectorspace, that is, we have
\[
S/I\iso \Dirsum_{i=1}^rS/P_i(-a_i).
\]
Each module $M_i$ in the filtration $\mathcal F$ is of the form $I_i/I$ where $I_i$ is a monomial ideal.  The monomials not belonging to $I_i$ form a $K$-basis
of $(S/I)/M_i=S/I_i$, and so $S/I=(S/I)/M_i\dirsum M_i$ decomposes naturally as a $K$-vectorspace. Identifying $S/P_i(-a_i)=M_i/M_{i-1}\subset S/I_{i-1}$ with
its image in $S$, we get $S/P_i(-a_i)=u_iK[Z_i]$ where $u_i=\prod_{j=1}^nx_j^{a_i(j)}$ and $Z_i=\{x_j\: j\not\in P_j\}$. Thus
\[
S/I= \Dirsum_{i=1}^ru_iK[Z_i].
\]
Any decomposition of $S/I$ as a direct sum of $K$-vectorspaces of the form $uK[Z]$ where $Z$ is a subset of $X=\{x_1,\ldots,x_n\}$ and $u$ is a monomial of
$K[X]$ is called a {\em Stanley decomposition}. Stanley decompositions have been studied in various combinatorial and algebraic contexts, see \cite{A},
\cite{HT}, and  \cite{MS}. Not all Stanley decompositions arise from prime filtrations, see \cite{MS}.

Stanley \cite{St1} conjectured that there always exists a Stanley decomposition  $S/I= \Dirsum_{i=1}^ru_iK[Z_i]$ such that $|Z_i|\geq \depth S/I$. In  \cite{A}
Apel studied cases in which Stanley's conjecture holds.

We conclude this section by showing

\begin{Theorem}
\label{stanley1} Let $I\subset S$  a monomial ideal, and suppose that $S/I$ is a multigraded pretty clean ring.  Then  Stanley's conjecture holds for $S/I$.
\end{Theorem}

\begin{proof}
Stanley's conjecture follows if we can show that there exist a multigraded prime filtration $\mathcal F$  of $S/I$ with factors $S/P_i(-a_i)$ such that $\depth
S/P_i\geq \depth S/I$.

Since $S/I$ is multigraded pretty clean, it follows from Corollary \ref{conclusion}  that all nonzero factors $D_i(S/I)/D_{i-1}(S/I)$ of the dimension
filtration are clean. Moreover, since $S/I$ is sequentially Cohen-Macaulay,  it follows from Proposition \ref{exti} that $\depth S/I=t$ where $t=\min\{i\:
D_i(S/I)/D_{i-1}(S/I)\neq 0\}$. Since $D_i(S/I)/D_{i-1}(S/I)$ is clean, we obtain a pretty clean filtration of $S/I$ as a refinement of the dimension
filtration by the clean filtrations of the factors $D_i(S/I)/D_{i-1}(S/I)$. Thus in this prime filtration each factor $S/P$ belongs to
$\Ass(D_i(S/I)/D_{i-1}(S/I))$ for some $i$. It follows that $\depth S/P\geq t$, as desired.
\end{proof}

\section{Prime filtrations and primary decompositions}

In this section we give another characterization of pretty clean modules in terms of primary decompositions.

\begin{Proposition}
\label{primary} Let $M$ be an $R$-module, and suppose $M$ admits the prime filtration
 ${\mathcal F}: (0)=M_0\subset M_1\subset \cdots \subset M_{r-1}\subset M_r=M$
 with $M_i/M_{i-1}\iso R/P_i$  for all $i$. Then for $j=1,\ldots, r$ there exist irreducible $P_j$-primary  submodules $N_j$ of $M$ such that
 $M_i=\Sect_{j=i}^rN_j$ for $i=0,\ldots, r$.
\end{Proposition}

In the proof of this result we shall need the following

\begin{Lemma}
\label{complement} Let $U\subset V\subset M$ be submodules of $M$ such that $V/U\iso R/P$ for some $P\in \Spec(R)$. Then there exists an irreducible submodule
$W$ of $M$ such that $U=V\sect W$.
\end{Lemma}

\begin{proof}
By Noetherian induction there exists a maximal submodule $W$ of $M$ such that
  $U=V\sect W$. We claim that $W$ is an irreducible submodule
  of $M$. Indeed, suppose that $W=W_1\sect W_2$. Then
  $U=(V\sect W_1)\sect(V\sect W_2)$
  is a decomposition of $U$ in $V$. However, $U$ is irreducible
 in $V$ since $V/U\iso R/P$. It follows that $V\sect W_1=U$ or $V\sect
 W_2=U$. Since $W$ was chosen to be maximal with this intersection property,
 we see that $W=W_1$ or $W=W_2$. Thus $W$ is irreducible, as desired.
\end{proof}

\begin{proof}[Proof of \ref{primary}] (a)\implies (b): Let ${\mathcal F}$  be a prime filtration as given in (a).
We show by decreasing induction on $i<r$  that for $j=i+1,\ldots, r$ there exist irreducible $P_j$-primary submodules $N_j$ of $M$ such that
$M_{i}=\Sect_{j={i+1}}^rN_j$.

For $i=r$ we may choose $N_r=M_{r-1}$, since $M/M_{r-1}\iso R/P_r$. Now let $1<i<r$, and assume that $M_{i}=\Sect_{j=i+1}^rN_j$ where $N_j$ is an irreducible
$P_j$-primary submodule of $M$ for $j=i+1,\cdots, r$. Since $M_{i}/M_{i-1}\iso R/P_{i}$, it follows by Lemma \ref{complement} that there exists an irreducible
submodule $N_{i}$ of $M$ such that $M_{i-1}=M_{i}\sect N_{i}$. Since $R/P_{i}\iso M_{i}/M_{i-1}=M_{i}/M_{i}\sect N_{i}\subset M/N_{i}$, it follows that
$\{P_{i}\}=\Ass(M_{i}/M_{i-1})\subset \Ass(M/N_{i})$. However $\Ass(M/N_{i})$ has only one element, therefore $\Ass(M/N_{i})=\{P_{i}\}$.
\end{proof}

\section{Clean filtrations and shellings}

In this section we recall the main result of the paper of Dress \cite{D} (see also \cite{Si}), and provide some extra information. Let $\Delta$ be a simplicial
complex on the vertex set $[n]=\{1,\ldots,n\}$. Recall that $\Delta$ is {\em shellable}, if the facets of $\Delta$ can be given a linear order $F_1,\ldots,
F_m$ such that for all $i,j$, $1\leq i<j\leq m$, there exists some $v\in F_i\setminus F_j$ and some $k<i$ with $F_i\setminus F_k=\{v\}$.

Note that we do {\em not} insist that $\Delta$ is pure, that is, that all facets of $\Delta$ have the same dimension. Sometimes such a shelling is called a
{\em non-pure shelling}.

Let $K$ be a field. The {\em Stanley-Reisner ring} of $K[\Delta]$ of $\Delta$ is the factor ring of $S=K[x_1,\ldots, x_n]$ modulo the ideal $I_\Delta$
generated by all squarefree monomials $x_{i_1}x_{i_2}\cdots x_{i_k}$ such that $\{i_1,\ldots, i_k\}$ is not a face of $\Delta$.

One has

\begin{Theorem}[Dress]
\label{dress} The simplicial complex $\Delta$ is shellable if and only if $K[\Delta]$ is a clean ring.
\end{Theorem}

For a subset of faces $G_1,\ldots, G_r$ of $\Delta$ we denote by $\langle G_1,\ldots, G_r\rangle$, the smallest subcomplex of $\Delta$ containing the faces
$G_1,\ldots, G_r$. With this notation, the shellability of $\Delta$ can also be characterized as follows: $\Delta$ is shellable if and only if the facets of
$\Delta$ can be ordered $F_1,\ldots, F_r$ such that for $i=2,\ldots,m$ the facets of $\langle F_1,\ldots, F_{i-1}\rangle\sect \langle F_i\rangle$
 are maximal proper faces of $\langle F_i\rangle$.

For $i\geq 2$ we denote by $a_i$ the number of facets of  $\langle F_1,\ldots, F_{i-1}\rangle\sect \langle F_i\rangle$, and set $a_1=0$. We call the
 $a_1,\ldots, a_r$ the sequence of {\em shelling numbers} of the given shelling of $\Delta$.

Set $P_{F_i}=(\{x_j\}_{j\not\in F_i})$. Then $I_\Delta=\Sect_{i=1}^rP_{F_i}$. Therefore, if $F_1,\ldots, F_r$ is a shelling of $\Delta$, then for
$i=2,\ldots,r$ we have
\[
\Sect_{j=1}^{i-1}P_{F_j}+P_{F_i} =P_{F_i}+(f_i).
\]
Here $f_i=\prod_k x_k$, where the product is taken over those $k\in F_i$ such that $F_i\setminus\{k\}$ is a facet of $\langle F_1,\ldots, F_{i-1}\rangle\sect
\langle F_i\rangle$. In particular it follows that  $\deg f_i$ equals   the $i$th shelling number $a_i$.

We obtain the following isomorphisms of graded $S$-modules
\begin{eqnarray*}
(\Sect_{j=1}^{i-1}P_{F_j})/(\Sect_{j=1}^{i}P_{F_j})&\iso &(\Sect_{j=1}^{i-1}P_{F_j}+P_{F_i})/P_{F_i}
=  (P_{F_i}+(f_i))/P_{F_i}\\
&\iso & (f_i)/(f_i)P_{F_i}\iso S/P_{F_i}(-a_i).
\end{eqnarray*}

The isomorphism $(P_{F_i}+(f_i))/P_{F_i}\iso (f_i)/(f_i)P_{F_i}$ results from the fact that $(f_i)\sect P_{F_i}=(f_i)P_{F_i}$ since the set of variables
dividing $f_i$ and the set of variables generating $P_{F_i}$ have no element in common. Thus we have shown

\begin{Proposition}
\label{shelling numbers} Let $\Delta$ be a shellable simplicial complex with shelling $F_1,\ldots, F_r$  and shelling numbers $a_1,\ldots, a_r$. Then
$(0)=M_0\subset M_1\subset\cdots\cdots M_{r-1}\subset M_r=K[\Delta]$ with
\[
M_i=\Sect_{j=1}^{r-i}P_{F_j}\quad \text{and}\quad M_i/M_{i-1}\iso S/P_{F_{r-i+1}}(-a_{r-i+1})
\]
is a clean filtration of $S/I_{\Delta}$.
\end{Proposition}

\section{Multicomplexes}

The aim of this and the next section is to extend the result of Dress to multicomplexes.  Stanley \cite{St} calls a subset $\Gamma\subset \NN^n$ a multicomplex
if for all $a\in \Gamma$ and all $b\in\NN^n$ with $b\leq a$, it follows that $b\in\Gamma$. The elements of $\Gamma$ are called {\em faces}.

What are the facets of $\Gamma$? We define on $\NN^n$ the partial order given by
\[ (a(1),\ldots, a(n))\leq (b(1),\ldots,b(n))\quad \text{if}\quad  a(i)\leq
b(i)\quad \text{for all}\quad i.
\]
An element $m\in \Gamma$ is called maximal if there exists no $a\in \Gamma$ with $a> m$. We denote by ${\mathcal M}(\Gamma)$ the set of maximal elements of
$\Gamma$. One would expect that ${\mathcal M}(\Gamma)$ is the set of facets of $\Gamma$. However  ${\mathcal M}(\Gamma)$ may be the empty set, for example for
$\Gamma=\NN^n$. To remedy this defect we will consider ``closed" subsets $\Gamma$ in $\NN^n_\infty$, where $\NN_\infty=\NN\union \{\infty\}$.

Let $a\in\Gamma$. Then $$\ip a=\{i\: a(i)=\infty\}$$ is called the {\em infinite part} of $a$. We first notice that

\begin{Lemma}
\label{finitem} Let $\Gamma\subset \NN^n_\infty$. Then ${\mathcal M}(\Gamma)$ is finite.
\end{Lemma}

\begin{proof}
Let $F\subset [n]$, and  set $\Gamma_F=\{a\in \Gamma\: \ip a =F\}.$ It is clear that if $a\in \Gamma_F$ is maximal in $\Gamma$ then $a$ is maximal in
$\Gamma_F$. Since there are only finitely many subsets $F$ of $[n]$, it suffices to show that $\Gamma_F$ has only finitely many maximal elements. Let
$[n]\setminus F=\{i_1,\ldots, i_k\}$ with $i_1<i_2<\cdots <i_k$. For each $a\in\Gamma_F$ we let $a'\in \NN^k$ be the integer vector with $a'(j)=a(i_j)$ for
$j=1,\ldots,k$. Now if $a$ and $b$ are two maximal elements in $\Gamma_F$ with $a\neq b$, then $a'$ and $b'$ are incomparable vectors, that is, $a'\not\leq b'$
and $b'\not\leq a'$. This implies that the set of monomials $\{x^{a'}\: a\in \Gamma_F,\; \text{$a$ maximal}\}$ is a minimal set of generators of the monomial
ideal they generate in $K[x_1,\ldots,x_k]$. Hence this set is finite. Thus  the set of maximal elements $\Gamma_F$ is finite for all $F\subset [n]$, and
${\mathcal M}(\Gamma)$ is finite.

\end{proof}

We say that a sequence of natural numbers  $a(i)$ has limit $\lim a(i)=\infty$, if for all integers $b$ there exists an integer $j$ such that $a(i)\geq b$ for
all $i\geq j$. Of course any non-decreasing sequence in $\NN$  has a limit -- either it is eventually constant, and this constant is its limit, or the limit is
$\infty$.

As usual we set $a\leq \infty$ for all $a\in \NN$. and  extend the partial order  on $\NN^n$ naturally to $\NN^n_\infty $.
 By what we just said it follows that any  sequence $a_i$, $i=1,2,\ldots$ of elements  in $\NN^n_\infty$ with $a_i\leq  a_{i+1}$ has a limit
-- the limit being taken componentwise.

\medskip
Let  $\Gamma\subset \NN_{\infty}^n$. The set  $\bar{\Gamma}$  of all $a\in\NN^n_\infty$ which are limits of ascending sequences in $\Gamma$ is called {\em the
closure} of $\Gamma§$.  It is clear that $\Gamma\subset \bar{\Gamma}$ and that $\bar{\bar{\Gamma}}=\bar{\Gamma}$.

\begin{Definition}
{\em A subset $\Gamma\subset \NN^n_\infty$ is called a {\em multicomplex} if
\begin{enumerate}
\item[(1)] for all $a\in\Gamma$ and all $b\in\NN^n_\infty$ with $b\leq a$ it follows that $b\in\Gamma$;

\item[(2)] $\Gamma=\bar{\Gamma}$.
\end{enumerate}}
\end{Definition}

The elements of a multicomplex are called {\em faces}. The next result shows that each face of a multicomplex is bounded by a face in ${\mathcal
M}(\Gamma)$.

\begin{Lemma}
\label{start} Let $\Gamma\subset \NN^n_\infty$ be a  set satisfying property $(1)$  of multicomplexes. Then the following conditions are equivalent:
\begin{enumerate}
\item[(a)] $\Gamma=\bar{\Gamma}$;

\item[(b)] for each $a\in\Gamma$ there exists $m\in{\mathcal M}(\Gamma)$ with $a\leq m$.
\end{enumerate}
\end{Lemma}

\begin{proof} (a)\implies (b): We proceed by induction on $n-|\ip a|$. If $n-|\ip a|=0$, then $a(i)=\infty$ for all $i$, and hence $a\in {\mathcal M}(\Gamma)$.
Suppose now that $n-|\ip a|>0$ and that there is no $m\in {\mathcal M}(\Gamma)$ with $a\leq m$.   Then there exists a strictly ascending sequence
$a=a_1<a_2<\ldots$ in $\Gamma$. Since $\Gamma=\Bar{\Gamma}$ it follows that  $b=\lim a_i\in \Gamma$. Obviously one has $n-|\ip b|<n-|\ip a|$. Hence by induction
hypothesis, there exists $m\in{\mathcal M}(\Gamma)$ with $b\leq m$, and thus  $a<m$.

(b)\implies (a): Let $a_i$, $i=1,2,\ldots$  be an ascending sequence in $\Gamma$. By assumption, there exist $m_i\in{\mathcal M}(\Gamma)$ with $a_i\leq m_i$.
Since ${\mathcal M}(\Gamma)$ is finite (see Lemma \ref{finitem}), there exists $i_0$ such that $m_i=m_{i_0}$ for all $i\geq i_0$. It follows that $a_i\leq
m_{i_0}$ for all $i$. Hence $\lim a_i\leq m_{i_0}$. In particular, $\lim a_i\in\Gamma$.
\end{proof}

Combining Lemma \ref{start} with Lemma \ref{finitem} we get

\begin{Corollary}
\label{m} Let $\Gamma\subset \NN^n_\infty$. Then $\Gamma$ is a multicomplex if and only if there exist finitely many elements $m_1,\ldots, m_r\in\NN^n_\infty$
such that
\[
\Gamma=\{a\in\NN^n_\infty\: a\leq m_i \text{ for some  }i=1,\ldots, r\}.
\]
\end{Corollary}

We have

\begin{Lemma}
\label{closure} Suppose  $\Gamma\subset \NN^n$ satisfies  property  $(1)$ of multicomplexes, then so does  $\bar{\Gamma}$.
\end{Lemma}

\begin{proof}
The statement is clear if $a\in\Gamma$. Suppose now that $a\in\bar{\Gamma}$, and let $a_i\in \Gamma$ be a non-descending sequence with $\lim a_i=a$. Let
$b_i(j)=\min\{a_i(j),b(j)\}$ for $j=1,\ldots,n$. Then  $b_i=(b_i(1),\ldots,b_i(n))\leq a_i$ for all $i$, and hence $b_i\in \Gamma$ for all $i$. Moreover,
$b=\lim b_i$ and so $b\in\bar{\Gamma}$.
\end{proof}

The lemma shows that if $\Gamma\subset \NN^n$ is a multicomplex in the sense of Stanley, then $\bar{\Gamma}\subset \NN^n_\infty$ is a multicomplex in our
sense. Moreover $\bar{\Gamma}\sect\NN^n=\Gamma$. Thus the assignment $\Gamma\mapsto \bar{\Gamma}$ establishes a bijection between these different concepts of
multicomplexes.

In the following we will use the term multicomplex only in our sense, that is,  we will always assume that $\Gamma=\bar{\Gamma}$.

\medskip
Note that $\Delta(\Gamma)=\{\ip a\: a\in\bar{\Gamma}\}$ is a simplicial complex on the vertex set $[n]=\{1,\ldots,n\}$. It is called the simplicial complex
associated to the multicomplex $\Gamma$.

The number  $\dim a=|\ip a|-1$ is called the {\em dimension of $a$}. The {\em dimension of $\Gamma$} is defined to be
\[
 \dim \Gamma= \max\{\dim a\: a\in \Gamma\}.
 \]
Obviously one has $\dim \Gamma =\dim \Delta(\Gamma)$.

\medskip
An element  $a\in \Gamma$ is called a {\em facet} of $\Gamma$  if for all $m\in {\mathcal M}(\Gamma)$ with  $a\leq m$ one has  $\ip a=\ip m$. The set of facets
of $\Gamma $ will be denoted by $\mathcal{F}(\Gamma)$. It is clear that $\mathcal{M}(\Gamma)\subset \mathcal{F}(\Gamma)$. The facets in $\mathcal{M}(\Gamma)$
are called {\em maximal facets}.

Consider for example the multicomplex $\Gamma\in\NN^2_\infty$ with faces $$\{a\: a\leq (0,\infty )\; \text{or}\; a\leq (2,0)\}.$$ Then ${\mathcal
M}(\Gamma)=\{(0, \infty), (2,0)\}$ and ${\mathcal F}(\Gamma)=\{(0,\infty), (2,0), (1,0)\}$. Besides its facets, $\Gamma$ admits the infinitely many faces
$(0,i)$ with $i\in \NN$.

\begin{Lemma}
\label{finite} Each multicomplex has a finite number of facets.
\end{Lemma}

\begin{proof}
Let $\Gamma$ be the given multicomplex.  Given $m\in{\mathcal M}(\Gamma)$. By \ref{finitem} it remains to show that the set
$$\{a\in \Gamma\: a\leq m\; \text{and}\; \ip a=\ip m\}$$ is
finite. But this is obviously the case since for each $i\not\in \ip m$ there are only $m(i)+1$ numbers $j\in \NN$ with $j\leq m(i)$.
\end{proof}

\begin{Lemma}
\label{intersection} An arbitrary intersection and a finite union of multicomplexes  is again a multicomplex.
\end{Lemma}

\begin{proof} Let $(\Gamma_i)_{i\in I}$ be a family of multicomplexes, and set $\Gamma=\Sect_{i\in I}\Gamma_i$.
If $a\in \Gamma$ and $b\leq a$, then obviously $b\in \Gamma$. Thus it  remains to show  that $\Gamma=\bar{\Gamma}$. Let $a_j$, $j=1,2,\ldots$ be an ascending
sequence in $\Gamma$. Since $\Gamma_i=\bar{\Gamma}_i$ for all $i\in I$, it follows that $\lim a_j\in \Gamma_i$ for all $i$, and hence $\lim a_i\in \Gamma$, as
desired.

On the other hand, suppose  $J=\{1,\ldots,k\}$ and let $\Gamma=\Union_{i=1}^k\Gamma_i$. Then $\Gamma$ satisfies obviously condition (1) of a multicomplex.

By Lemma \ref{finitem} the sets  ${\mathcal M}(\Gamma_i)$ are finite, and  $\Union_{i=1}^k\Gamma_i$ is the set of all $a\in\NN^n_\infty$ for which there exists
$j\in J$ and $m\in {\mathcal M}(\Gamma_j)$ such that $a\leq m$. Thus it follows from Corollary \ref{m} the $\Gamma$ is a multicomplex.
\end{proof}

\begin{Corollary}
\label{closure} Let $A\subset \NN^n_\infty$ be an arbitrary subset of $\NN^n_\infty$. Then there exists a unique smallest multicomplex $\Gamma(A)$ containing
$A$.
\end{Corollary}

Let $\Gamma$ be a multicomplex, and  let $I(\Gamma)$ be the $K$-subspace in $S=K[x_1,\ldots,x_n]$ spanned by all monomials $x^a$ such that $a\not\in \Gamma$.
Note that if $a\in \NN^n$ and $b\in \NN^n\setminus\Gamma$, then $a+b\in \NN^n\setminus \Gamma$, that is, if $x^a\in I(\Gamma)$ then $x^ax^b\in I(\Gamma)$ for
all $x^b\in S$. In other words, $I(\Gamma)$ is a monomial ideal. In particular, the monomials $x^a$ with $a\in \Gamma$ form a $K$-basis of $S/I(\Gamma)$.

For example for the above multicomplex $\Gamma=\{a\: a\leq (0,\infty )\; \text{or}\; a\leq (2,0)\}$ in $\NN^2_\infty$ we have $I(\Gamma)=(x_1^3,x_1x_2)$.

Conversely, given an arbitrary monomial ideal $I\subset S$, there is a unique multicomplex $\Gamma$ with $I=I(\Gamma)$. Indeed, let $A =\{a\in \NN^n\:
x^a\not\in I\}$; then  $\Gamma=\Gamma(A)$.

\medskip
The monomial ideal of a multicomplex behaves with respect to intersections and unions of multicomplexes as follows:

\begin{Lemma}
\label{ideal} Let $\Gamma_j$, $j\in J$ be a family of multicomplexes. Then
\begin{enumerate}
\item[(a)] $I(\Sect_{j\in J}\Gamma_j)=\sum_{j\in J} I(\Gamma_j)$, \item[(b)] if $J$ is finite, then $I(\Union_{j\in J}\Gamma_j)=\Sect_{j\in J}I(\Gamma_j)$.
\end{enumerate}
\end{Lemma}
\noindent
 Next we describe the relationship between simplicial
complexes and multicomplexes. Let $\Delta$ be a simplicial complex on the vertex set $[n]$. To each facet $F\in \Delta$ we associate the element
$a_F\in\NN^n_\infty$ with
\[
a_F(i)=\left\{ \begin{array}{lll} \infty, & \text{if} & i\in F\,\\
0, & \mbox{if} & i\not\in F, \end{array} \right.
\]
Then $\{a_F\: F\in \Delta\}$ is the set of facets of a multicomplex $\Gamma(\Delta)$, and $I(\Gamma(\Delta))=I_\Delta$, where $I_\Delta$ is the Stanley-Reisner
ideal of $\Delta$. Moreover one has $\dim \Gamma=\dim \Delta(\Gamma)$.

\medskip
For a multicomplex $\Gamma$ and $a\in \Gamma$ we  let $P_a$ be the prime ideal generated by all $x_i$ with $i\not\in \ip a$. Thus $P_a$ is generated by all
$x_i$ with $a(i)\in \NN$.

\begin{Lemma}
\label{irred} Let $\Gamma$ be a multicomplex. The following statements are equivalent:
\begin{enumerate}
\item[(a)] $\Gamma$ has just one maximal facet $a$;

\item[(b)]  $I(\Gamma)$ is an irreducible ideal.
\end{enumerate}
If the equivalent conditions hold, then $I(\Gamma)$ is generated by $\{x_i^{a(i)+1}\: i\in [n]\setminus \ip a\}$. In particular, $I(\Gamma)$ is a
$P_a$-primary ideal.
\end{Lemma}

\begin{proof}
If $a$ is the unique maximal facet of $\Gamma$ then
\[
I(\Gamma)=(x^b\: b\in\NN^n,\; b(i)>a(i)\; \text{for some $i$}) =(x_i^{a(i)+1}\: i\in [n]\setminus \ip a).
\]
Conversely, if $I(\Gamma)$ is irreducible, then according to  \cite[Theorem 5.1.16]{Vi} there exists a  subset $A\subset \{1,\ldots,n\}$  and  for each  $i\in
A$ an integer  $a_i>0$  such that $I(\Gamma)=(x_i^{a_i}:i\in A, a_i>0)$. Set $a(i)=a_i-1$ for $i\in A$ and $a(i)=\infty$ for $i\not \in A$. Then $a$ is the
unique facet of $\Gamma$.
\end{proof}

\begin{Corollary}
\label{irrprime} Let $\Gamma\subset \NN_{\infty}^n$ be a multicomplex with just one facet $a$. Then $I(\Gamma)=P_a$.
\end{Corollary}

\begin{proof}
Suppose $a(i)\neq 0$ for some $i\not\in  \ip a$. Then $a-e_i$ is a facet, different from $a$. Here $e_i$ is the canonical $i$th unique vector. Thus we see that
$a(i)\in \{0,\infty\}$ for $i=1,\ldots,n$, so that $I(\Gamma)=I(\Gamma(a))=P_a$.
\end{proof}

The next result describes how the maximal facets of a multicomplex $\Gamma$ are related to the irreducible components of $I(\Gamma)$.

\begin{Proposition}
\label{irrdec} Let $\Gamma\subset \NN_{\infty}^n$ be a multicomplex, and $a_1,\ldots,a_r$ its maximal facets. Then
$$I(\Gamma)=\Sect_{j=1}^rI(\Gamma(a_j))$$
 is the unique irredundant
irreducible decomposition of $I(\Gamma)$ in $S=K[x_1,\ldots,x_n]$.

Conversely, let $I\subset S$ be a monomial ideal, $I=\Sect_{j=1}^rI_j$  the unique irredundant irreducible decomposition of $I$ in $S$, and let $\Gamma$ be the
multicomplex with $I(\Gamma)=I$. Then $\Gamma$ has $r$ maximal facets $a_1,\ldots, a_r$ which can be labelled such that
$$I(\Gamma(a_j)) =I_j\quad \text{for}\quad j=1,\ldots, r.$$
\end{Proposition}

\begin{proof}
Since $\Gamma=\Union_{i=1}^r\Gamma(a_i)$, it follows from Lemma \ref{ideal} that $I(\Gamma)=\Sect_{j=1}^rI(\Gamma (a_j))$. That each $I(\Gamma(a_i))$ is
irreducible, we have seen in Lemma \ref{irred}.

Conversely, let $I=\Sect_{j=1}^rI_j$ be the unique irredundant irreducible decomposition of $I$, and let $\Gamma_j$ be the unique multicomplex with
$I(\Gamma_j)=I_j$. By Lemma \ref{irred}, each $\Gamma_j$ has exactly one maximal facet, say $a_j$. Hence $\Gamma_j=\Gamma(a_j)$ for $j=1,\ldots, r$.

Let $\Gamma$ be the unique multicomplex with $I(\Gamma)=I$. Then since $I(\Gamma)=\Sect_{i=1}^rI((\Gamma(a_j))$, it follows from Lemma \ref{ideal} that
$I(\Gamma)=I(\Gamma( a_1,\ldots,a_r))$, and hence that $\Gamma=\Gamma(a_1,\ldots,a_r)$. Each of the $a_j$ is a maximal  facet of $\Gamma$, because if there
would be an inclusion among them, then there would also be an inclusion among the $I_j$, contradicting the minimality of the decomposition.
\end{proof}

\begin{Corollary}
\label{dimension} Let $\Gamma$ be a multicomplex. Then $\dim S/I(\Gamma)=\dim \Gamma+1$.
\end{Corollary}

\begin{proof}
By the preceding proposition it suffices to prove the assertion in case that $\Gamma$ has just one maximal facet, say $a$. Suppose that $\dim \Gamma=d-1$. We
may, then assume that $a(i)=\infty$ for $i\geq n-d+1$. Then $I(\Gamma)=(x_1^{a(1)+1},\ldots,x_{n-d}^{a(n-d)+1})$, and $\dim S/I(\Gamma)=d$.
\end{proof}

\medskip
\noindent Finally  we will show that the facets of a multicomplex $\Gamma$ correspond to the standard pairs of $I=I(\Gamma)$ introduced by Sturmfels, Trung and
Vogel \cite{STV}: let $u$ be a monomial of $S=K[x_1,\ldots, x_n]$. Then we set $\supp(u)=\{x_i\: x_i\text{ divides } u\}$. A pair $(u,Z)$ where $u$ is a
monomial and $Z$ is a subset of the set of variables $X=\{x_1,\ldots, x_n\}$ is called {\em admissible} if no $x_i\in Z$ divides $u$, that is, if
$\supp(u)\sect Z=\emptyset$. The set of admissible pairs is partially ordered as follows:
\[
(u,Z)\leq(u',Z')\quad\iff\quad \text{$u$ divides $u'$}\quad \text{and}\quad \supp(u'/u)\union Z'\subset Z.
\]
An admissible pair $(u,Z)$ is called {\em standard} with respect to $I$, if $u K[Z]\sect I=\{0\}$, and $(u, Z)$ is minimal with this property. The set of
standard pairs with respect to $I$ is denoted by $\std(I)$.

\medskip
\noindent For a monomial $u\in S$, with $u=\prod_{i=1}^nx_i^{a_i}$ we set $\log u= (a_1,\ldots, a_n)$, and for a subset $Z\subset X$ we let $c(Z)\in
\NN_\infty^n$ the element with
\[
c(Z)(i) = \left\{ \begin{array}{lll} \infty, & \text{if} & x_i\in Z,\\
0, & \mbox{if} & x_i\not\in Z. \end{array} \right.
\]

With this notation we have

\begin{Lemma}
\label{pairs} Let $I\subset S$ be a monomial ideal, and $\Gamma$ the multicomplex associated with $I$. Then the standard pairs with respect to $I$ correspond
bijectively to the facets of $\Gamma$. The bijection is established by the following assignment:
\[
\std(I)\To {\mathcal F}(\Gamma),  \quad (u,Z)\mapsto \log u+c(Z).
\]
\end{Lemma}

\begin{proof} Let $\mathcal A$ be the set of admissible pairs. Since $\supp u\sect Z=\emptyset$ for  $(u,Z)\in \mathcal A$ it follows that the map
\[
{\mathcal A}\To \NN^n_\infty,\quad (u,Z)\mapsto \log u+c(Z)
\]
is injective. Moreover, for each $(u,Z)\in \mathcal A$ we have
\[
uK[Z]\sect I=\{0\}\iff \log u+c(Z)\in \Gamma.
\]
Now let $(u,Z)\in\std(I)$, and set  $a=\log u+c(Z)$. Let $m\in {\mathcal M}(\Gamma)$ with $a\leq m$. Suppose that $\ip a \neq \ip m$. Then there exists $i$
such that $a(i)<m(i)=\infty$. Let $v=u/x_i^{a(i)}$ and $W=Z\union \{x_i\}$. Then  $(v,W)<(u,Z)$ and $v\cdot K[W]\sect I=\{0\}$, a contradiction. Therefore,
$a\in{\mathcal F}(\Gamma)$.

Conversely let $a\in {\mathcal F}(\Gamma)$. Set $u=\prod_{i\not\in\ip(a)}x_i^{a(i)}$ and $Z=\{x_i\: i\in\ip a\}$.  Then $(u,Z)\in\mathcal A$ and $a=\log
u+c(Z)$. Since $a\in\Gamma$ it follows that $u\cdot K[Z]\sect I=\{0\}$. Suppose that $(u,Z)$  is not minimal with this property. Then there exists $(v,W)\in
\mathcal A$ with  $v\cdot K[W]\sect I=\{\ 0\}$ and $(v,W)<(u,Z)$, and we have

\begin{enumerate}
\item[(1)] $b=\log v+c(W)\in \Gamma$;

\item[(2)] $v$ divides $u$;

\item[(3)] $\supp(u/v)\union Z\subset W$.
\end{enumerate}

The properties (2) and (3) imply that $a(i)=b(i)$ for all $i$ such that $b(i)<\infty$. Thus $a\leq b$, and $a=b$ if and only if $\ip a=\ip b$. However since
$a\neq b$, we have $\ip a\neq \ip b$. By property (1) there exists  $m\in {\mathcal M}(\Gamma)$ with $b\leq m$. Then $a\leq m$ and $\ip b\subset \ip m$. In
particular, $\ip a \neq \ip m$. It follows that $a\not \in {\mathcal F}(\Gamma)$, a contradiction.
\end{proof}

\section{Pretty clean filtrations and shellable multicomplexes}

In this section we introduce shellable multicomplexes and show how this concept is related to clean filtrations. Our concept of shellability is a translation
of Corollary  \ref{primary1}  into the language of multicomplexes. In that corollary we characterized pretty clean filtrations in terms of primary
decompositions. Here we need a refined multigraded  version of this result.

\begin{Proposition}
\label{multiprimary} Let $S=K[x_1,\ldots, x_n]$ be  the polynomial ring, and $I\subset S$ a monomial ideal. The following conditions are equivalent:
\begin{enumerate}
\item[(a)]  $S/I$ admits a multigraded prime filtration
 ${\mathcal F}: (0)=M_0\subset M_1\subset \cdots \subset M_{r-1}\subset M_r=S/I$
 such that $M_i/M_{i-1}\iso S/P_i(-a_i)$  for all $i$;
\item[(b)] there exists a chain of monomial ideals $I=I_0\subset I_1\subset \cdots \subset I_r=S$ and monomials $u_i$ of multidegree $a_i$ such that
$I_i=(I_{i-1},u_i)$ and $I_{i-1}:u_i=P_i$;
\end{enumerate}
If the equivalent conditions hold, then there exist irreducible monomial ideals $J_1,\ldots J_r$ such that $I_i=\Sect_{j=i+1}^r J_j$ for $i=0,\ldots, r$.
Moreover, if the prime filtration is pretty clean, then this set of irreducible ideals $\{J_1,\ldots, J_r\}$  is uniquely determined. In fact, this set
corresponds bijectively to the set of facets of the multicomplex associated with $I$.
\end{Proposition}

\begin{proof} The statements (a) and (b) are obviously equivalent, while the existence of of the irreducible ideals $J_i$  is just the multigraded version of Proposition
\ref{primary}.

Now we assume that the prime filtration $\mathcal F$ is pretty clean. Since $J_i$ is an irreducible monomial ideal,  it follows that $J_i=\Gamma(a_i)$ for some
$a_i\in\NN^n_\infty$, see Lemma \ref{irred}. We claim that ${\mathcal A}=\{a_1,\ldots, a_r\}$ is the set of facets of  the unique multicomplex $\Gamma$ with
$I=I(\Gamma)$.

We first show that all $a_j$ are facets of $\Gamma$. Note that ${\mathcal M}(\Gamma)\subset \mathcal A$. Indeed,  by Proposition \ref{irrdec} we have that
$$I(\Gamma)=\Sect _{a\in {\mathcal M}(\Gamma)}I(\Gamma(a))$$
is the unique irredundant decomposition of $I(\Gamma)$ into  irreducible ideals. Since from any redundant such decomposition, like the decomposition
$I=\Sect_{j=1}^rJ_j$, we obtain an irredundant by omitting redundant components we obtain the desired inclusion.

We also see that for each $J_j$ there exists a maximal facet $a$ of $\Gamma$ such that $I(\Gamma(a))\subset J_j$, that is, for each  $a_j\in \mathcal A$ there
exists a maximal facet $a$ of $\Gamma$ such that $a_j\leq a$. We claim that $\ip a_j=\ip a$, in other words, that $P_a= P_j$. In fact, since $a\in \mathcal A$
as we have just seen, there exists an integer $i$ such that $a=a_i$, and hence $I(\Gamma(a))=J_i$ is $P_i$-primary, and $P_i\subset P_j$. Suppose that $P_i\neq
P_j$. Then, since $\mathcal F$ is pretty clean, we conclude that $i>j$. It follows that $\Sect_{t>j}J_t=\Sect_{t\geq j}J_t$, contradicting (b).

Thus we have shown that all elements of $\mathcal A$ are facets of $\Gamma$. Next we prove that $r=|{\mathcal F}(\Gamma)|$. This then implies that ${\mathcal
A}={\mathcal F}(\Gamma)$, and that the elements of $\mathcal A$ are pairwise distinct.

We know from Corollary \ref{independent} that $r$ equals  the arithmetic degree of $S/I$. On the other hand we have shown in Lemma \ref{pairs} that the facets
of $\Gamma$ correspond to the standard pairs of $I$. In \cite[Lemma 3.3]{STV} it is shown that the number of standard pairs of $I$ is equal to  the arithmetic
degree of $S/I$ as well. Thus $|{\mathcal F}(\Gamma)|=r$, as desired.
\end{proof}

\medskip
\noindent
 In Section 6 we have considered the Stanley decomposition of $S/I$ into subspaces of the form $uK[Z]$ where $u$ is a monomial in the variables
$X=\{x_1,\ldots, x_n\}$ and $Z\subset X$. We call $S\subset \NN^n_\infty$ a {\em Stanley set} if there exists $a\in \NN^n$ and $m\in \NN^n_\infty$ with
$m(i)\in\{0,\infty\}$ such that $S=a+S^*$, where  $S^*=\Gamma(m)$ . The {\em dimension of $S$} is defined to be $\dim\Gamma(m)$. Obviously Stanley sets
correspond to subspaces of the form $uK[Z]$.

\begin{Definition}
{\em A multicomplex $\Gamma$ is {\em shellable} if the facets of $\Gamma$ can be ordered $a_1,\ldots, a_r$ such that
\begin{enumerate}
\item[(1)] $S_i= \Gamma(a_i)\setminus\Gamma(a_1,\ldots, a_{i-1})$ is a Stanley set for $i=1,\ldots,r$, and

\item[(2)] whenever $S_i^*\subset S_j^*$, then $S_i^*=S_j^*$ or $i>j$.
\end{enumerate}
Any order of the facets satisfying (1) and (2) is called a {\em shelling} of $\Gamma$}
\end{Definition}

The next result   shows that our definition of shellability of multicomplexes extends the classical concept of shellability of simplicial complexes.

\begin{Proposition}
\label{extend} Let $\Delta$ be a simplicial complex with facets $F_1,\ldots, F_r$, and $\Gamma$ be the multicomplex   with facets $a_{F_1},\ldots, a_{F_m}$.
Then $F_1,\ldots, F_m$ is a shelling of $\Delta$ if and only if $a_{F_1},\ldots, a_{F_r}$ is a shelling of $\Gamma$.
\end{Proposition}

\begin{proof}
We denote by $e_i$ the $i$th standard unit vector in $\NN^n$, and set $\Gamma_i=\Gamma(a_{F_i})$. Then

\begin{eqnarray*}
\Gamma(a_{F_i})\setminus\Gamma(a_{F_1},\ldots,a_{F_{i-1}})&=&\Sect_{j=1}^{i-1}(\Gamma(a_{F_i})\setminus\Gamma(a_{F_j}))\\
&=&\Sect_{j=1}^{i-1}(\Union_{k\in F_i\setminus F_j}(e_k+\Gamma_i))
\end{eqnarray*}
We notice that
\[
(e_k+\Gamma_i)\sect (e_l+\Gamma_i)=\left\{ \begin{array}{lll} e_k+\Gamma_i, & \text{if} & k=l,\\
e_k+e_l+\Gamma_i, & \mbox{if} & k\neq l, \end{array} \right.
\]
Thus
\[
\Gamma(a_{F_i})\setminus\Gamma(a_{F_1},\ldots,a_{F_{i-1}})=\Union_{L\in\mathcal L}(e_L+\Gamma_i),
\]
where  \[ {\mathcal L}=\{\{k_1,\ldots, k_{i-1}\}\: k_j\in F_i\setminus F_j \text{ for } j=1,\ldots,{i-1}\}
\]
and where $e_L=\sum_{j\in L}e_{j}$ for each $L\in\mathcal L$.

The union
\[
\Union_{L\in\mathcal L}(e_L+\Gamma_i)
\]
is a Stanley set if and only if there exists $L\in\mathcal L$ such that $e_{L'}+\Gamma_i\subset e_L+\Gamma_i$ for  all $L'\in \mathcal L$, and this is the case
if and only if there exists $L\in\mathcal L$ such that $L\subset L'$ for all $L'\in \mathcal L$.

We claim that the last condition is equivalent to the condition that all facets of $\langle F_i\rangle\sect \langle F_1,\ldots, F_{i-1}\rangle$ are maximal
proper subfaces of $\langle F_i\rangle$.

Suppose first that there is a set $L_0\in\mathcal L$ which is minimal under inclusion. We may assume that $L_0=[m]$. Let $k\in [m]$ and assume that all sets
$F_i\setminus F_j$  which contain $k$ have more than one element. Then for each such set we can pick $k_j\in F_i\setminus F_j$ with $k_j\neq k$, and hence
there exists $L\in\mathcal L$ which does not contain $k$, a contradiction, since $k\in L_0\subset L$. Thus for each $k\in L_0$ there exists an integer $j_k\in
[i-1]$ such that $F_i\setminus F_{j_k}=\{k\}$. Now let $j\in [i-1]$ be arbitrary. If $|F_i\setminus F_j|=1$, then by definition of the sets $L$, the set
$F_i\setminus F_j$ is a subset of each $L$, and in particular of $L_0$. Thus we see that the subfaces of $\langle F_i\rangle\sect \langle F_1,\ldots,
F_{i-1}\rangle$ of codimension $1$ are exactly the faces $F_i\setminus\{k\}$ for $k=1,\ldots, m$. Suppose now there exists  $j\in [i-1]$ for which $F_i\sect
F_j$ is not contained in any of these codimension 1 subfaces of $F_i$ (in which case not all facets of $\langle F_i\rangle\sect \langle F_1,\ldots,
F_{i-1}\rangle$ would be  maximal proper subfaces of $\langle F_i\rangle$.). Then $k\not\in F_i\setminus F_j$ for $k=1,\ldots, m$, and hence $(F_i\setminus
F_j)\sect L_0=\emptyset$. This a  contradiction, since any $L\subset \mathcal L$ contains an element of $F_i\setminus F_j$.

Conversely, suppose that  all facets of $\langle F_i\rangle\sect \langle F_1,\ldots, F_{i-1}\rangle$ are maximal proper subfaces of $\langle F_i\rangle$.
Then there exist $j_1,\ldots, j_m\in [i-1]$ such that $|F_i\setminus F_{j_k}|=1$, and for any $j\in [i-1]$ there exists $k\in [m]$ such that $F_i\setminus
F_{j_k}\subset F_i\setminus F_j$. For simplicity we may assume that $F_i\setminus F_{j_k}=\{k\}$ for $k=1,\ldots,m$. Then obviously  $L_0\in\mathcal L$ and
$L_0\subset L$ for any other $L\in\mathcal L$.
\end{proof}

\begin{Remark}
{\em  Condition (2)  in the definition of shellability  is superfluous  in case $\Gamma$ is the multicomplex  corresponding to a simplicial complex, because in
this case the sets $S_i^*$ correspond to the minimal prime ideals of $I(\Gamma)$, and hence there is no inclusion among them.}
\end{Remark}

\medskip
As an extension of the theorem of Dress we now show
\begin{Theorem}
\label{multi2} The multicomplex $\Gamma$ is shellable if and only if $S/I(\Gamma)$ is a multigraded pretty clean ring.
\end{Theorem}

\begin{proof}
Let $a_1,\ldots, a_r$ be the facets of $\Gamma$, and let $J_{j}=\Gamma(a_j)$ for $j=1,\ldots,r$. Then $J_j$ is an irreducible monomial ideal, and
$I(\Gamma(a_1,\ldots, a_{i}))=\Sect_{j=1}^{i}J_j$. We set $I_i=\Sect_{j=1}^{r-i+1}J_j$ and $M_i=I_i/I$ for $i=0,\ldots, r$, $I=I(\Gamma)$. Then ${\mathcal F}\: (0)=M_0\subset
M_1\subset \cdots\subset M_r=S/I$ is a multigraded filtration of $S/I$.

Since $\Gamma(a_i)\setminus \Gamma(a_1,\ldots, a_{i-1})=\Gamma(a_1,\ldots, a_{i})\setminus \Gamma(a_1,\ldots, a_{i-1})$, we see that $b\in \Gamma(a_i)\setminus
\Gamma(a_1,\ldots, a_{i-1})$ if and only $x^b\in \Sect_{j=1}^{i-1}iJ_j\setminus \Sect_{j=1}^{i}J_j$. In other words, the monomials $x^b$ with $b\in
\Gamma(a_i)\setminus \Gamma(a_1,\ldots, a_{i-1})$ form $K$-basis of the factor module $I_i/I_{i-1}=M_{i}/M_{i-1}$.

The discussion at the end of Section 6 shows that $\mathcal F$ is a  prime filtration if and only if $M_i/M_{i-1}$ as monomial vectorspace is isomorphic to
$uK[Z]$ for some monomial $u\in S$ and some subset $Z\subset \{x_1,\ldots,x_n\}$. Consequently, $\mathcal F$ is a prime filtration if and only
$\Gamma(a_i)\setminus \Gamma(a_1,\ldots, a_{i-1})$ is a Stanley set for all $i=1,\ldots,r$. Hence the theorem follows from Proposition \ref{multiprimary}.
\end{proof}

Let $K$ be field, and let $S=K[x_1,\ldots,x_n]$ be the polynomial ring. We call a multicomplex $\Gamma\subset \NN^n_\infty$ {\em Cohen-Macaulay} or {\em
sequentially Cohen-Macaulay over $K$} if $S/I(\Gamma)$ has the corresponding property.

$\Gamma$ is simply called Cohen-Macaulay, or sequentially Cohen-Macaulay, if $S/I(\Gamma)$ has the corresponding property over any field.

\begin{Corollary}
\label{properties} Let $\Gamma$ be a shellable multicomplex. Then $\Gamma$ is sequentially Cohen-Macaulay. If moreover, all facets of $\Gamma$ have the same
dimension, then $\Gamma$ is Cohen-Macaulay.
\end{Corollary}

\begin{proof}
Theorem \ref{multi2} implies that $S/I(\Gamma)$ is pretty clean. Hence the assertions follow from Theorem \ref{sequentially}.
\end{proof}

\begin{Corollary}
\label{criterion} A multicomplex $\Gamma$ is shellable if and only if there exists an order $a_1,\ldots, a_r$ of the facets such that for $i=1,\ldots,r$ the
sets $S_i=\Gamma(a_i)\setminus\Gamma(a_1,\ldots, a_{i-1})$ are Stanley sets with $\dim S_1\geq \dim S_2\geq \ldots\geq \dim S_r$.
\end{Corollary}

\begin{proof} Suppose the conditions of the corollary are satisfied, and that $S_i^*\subset S_j^*$ for some $i<j$. Then,  since $\dim S_i
\geq S_j$, it follows that $S_i^*=S_j^*$. Thus $\Gamma$ is shellable.

Conversely, suppose hat $\Gamma$ is shellable. Then $S/I$ is pretty clean. Thus by Corollary \ref{conclusion} the non-zero factors of the dimension filtration
are clean. Refining the dimension filtration by the clean filtrations of the factors we obtain a pretty clean filtration with $\dim S_1\geq \dim S_2\geq
\ldots\geq \dim S_r$.
\end{proof}

\end{document}